\documentclass[11pt]{article}

\usepackage{amssymb}

\makeatletter
\@addtoreset{equation}{section}

\newtheorem{theorem}{Theorem}[section]
\newtheorem{proposition}[theorem]{Proposition}
\newtheorem{corollary}[theorem]{Corollary}
\newtheorem{lemma}[theorem]{Lemma}
\newtheorem{definition}[theorem]{Definition}
\newtheorem{remark}[theorem]{Remark}

\newtheorem{example}[theorem]{Example}

\newcommand{\tends}[1]{{\displaystyle\mathop{\to}_{#1}}}

\def\endverif{\nopagebreak\newline\mbox{\ }\hfill\rule{2mm}{2mm}}
\def\endveriff{\nopagebreak\mbox{\ }\hfill\rule{2mm}{2mm}}

\def\1{{\bf 1}}

\def\2{{1/2}}

\def\C{{\mathbb C}}
\def\N{{\mathbb N}}
\def\R{{\mathbb R}}

\def\Z{{\mathbb Z}}

\def\A{{\cal A}}

\def\H{{\cal H}}

\def\RR{{\cal R}}

\def\eps{\varepsilon}

\def\Ad{{\rm Ad}\,}
\def\CPA{{\rm CPA}}

\def\Ker{{\rm Ker}}

\def\rank{{\rm rank}\,}

\def\span{{\rm span}}

\def\supp{{\rm supp}\,}
\def\Tr{{\rm Tr}}

\begin{document}

\title{\bf The variational principle for a class of asymptotically
abelian C$^*$-algebras}

\author{Sergey Neshveyev\footnote{Partially supported by NATO grant
SA (PST.CLG.976206)5273}\\
{\it B. Verkin Institute for Low Temperature Physics and Engineering}\\
{\it 47, Lenin Ave., 310164, Kharkov, Ukraine}\\
\\
and\\
\\
Erling St{\o}rmer\footnote{Partially supported by the Norwegian
Research Council}\\
{\it Department of Mathematics, University of Oslo,}\\
{\it P.O. Box 1053, Blindern, 0316 Oslo, Norway}}

\date{}

\maketitle

\begin{abstract}
Let $(A,\alpha)$ be a C$^*$-dynamical system. We introduce the notion
of pressure $P_\alpha(H)$ of the automorphism $\alpha$ at a
self-adjoint operator $H\in A$. Then we consider the class of AF-systems
satisfying the following condition: there exists a dense
$\alpha$-invariant $*$-subalgebra $\A$ of~$A$ such that for all
pairs $a,b\in\A$ the C$^*$-algebra they generate is finite
dimensional, and there is $p=p(a,b)\in\N$ such that
$[\alpha^j(a),b]=0$ for $|j|\ge p$. For systems in this class we
prove the variational principle, i.e. show that $P_\alpha(H)$ is the
supremum of the quantities $h_\phi(\alpha)-\phi(H)$, where $h_\phi(\alpha)$
is the Connes-Narnhofer-Thirring dynamical entropy of $\alpha$ with
respect to the $\alpha$-invariant state~$\phi$. If $H\in\A$, and
$P_\alpha(H)$ is finite, we show that any state on which the supremum
is attained is a KMS-state with respect to a one-parameter
automorphism group naturally associated with~$H$. In particular,
Voiculescu's topological entropy is equal to the supremum of
$h_\phi(\alpha)$, and any state of finite maximal entropy is a trace.
\end{abstract}

\section{Introduction} \label{1}

The variational principle has over the years attracted much attention
both in classical ergodic theory, see e.g. \cite{W}, and in the
C$^*$-algebra setting of quantum statistical mechanics, see e.g.
\cite{BR2}. In the years around 1970 there was much progress in the
study of spin lattice systems. In that case one is given for each
point $x\in\Z^\nu$ ($\nu\in\N$) the algebra of all linear operators
$B(\H)_x$ on a finite dimensional Hilbert space $\H$, and if
$\Lambda\subset\Z^\nu$ one defines the C$^*$-algebra
$A(\Lambda)=\otimes_{x\in\Lambda}B(\H)_x$. One then considers the
UHF-algebra $A=\overline{\cup_{\Lambda\subset\Z^\nu}A(\Lambda)}$ with space
translations $\alpha$, and studies mean entropy $s(\phi)$ of invariant states
and the corresponding variational principle
\begin{equation} \label{e1.1}
P(H)=\sup_\phi(s(\phi)-\phi(H))
\end{equation}
together with the KMS-states defined by a natural derivation
associated with a self-adjoint operator $H$, see \cite[Chapter~6]{BR2}.

With the development of dynamical entropy of automorphisms of
C$^*$-algebras \cite{CS, CNT, V} it was natural to replace the mean
entropy $s(\phi)$ by dynamical entropy $h_\phi(\alpha)$. This was
done by Narnhofer \cite{N},
who considered KMS-properties of the states on which the quantity
$h_\phi(\alpha)-\phi(H)$ attains its maximal value.
Then Moriya
\cite{M} showed that one can replace $s(\phi)$ by the CNT-entropy
$h_\phi(\alpha)$ in (\ref{e1.1}) and get the same result, i.e.
\begin{equation} \label{e1.2}
P(H)=\sup_\phi(h_\phi(\alpha)-\phi(H)).
\end{equation}

If one wants to study the variational principle and equilibrium
states for more general C$^*$-dynamical systems, the mean entropy is
usually not well defined, and it is necessary to consider dynamical
entropy. However, in order to define time
translations and extend the theory of spin lattice systems rather strong
assumptions of asymptotic abelianness are needed.

In the present paper we shall study a restricted class of
asymptotically abelian systems $(A,\alpha)$, namely we shall assume
that $A$ is a unital separable C$^*$-algebra, and $(A,\alpha)$ is
asymptotically abelian with locality, i.e. there exists a dense
$\alpha$-invariant $*$-subalgebra $\A$ of~$A$ such that for all
pairs $a,b\in\A$ the C$^*$-algebra they generate is finite
dimensional, and there is $p=p(a,b)\in\N$ such that
$[\alpha^j(a),b]=0$ for $|j|\ge p$. In particular, $A$ is an
AF-algebra. Examples of such systems are described in~\cite{S}. They are
all different shifts, on infinite tensor products of the same AF-algebra
with itself, on the Temperley-Lieb algebras, on towers of relative
commutants, on binary shifts algebras
defined by finite subsets of the natural numbers.

In Section \ref{2} we define the pressure $P_\alpha(H)$ of~$\alpha$
with respect to a self-adjoint operator $H\in A$. This can be done in
any unital C$^*$-dynamical system $(A,\alpha)$ with $A$ a nuclear
C$^*$-algebra, and follows closely Voiculescu's definition of
topological entropy $ht(\alpha)$ in \cite{V}. The main difference is
that he considered $\rank B$ of a finite dimensional C$^*$-algebra
$B$, while we consider quantities of the form $\Tr_B(e^{-K})$ for $K$
self-adjoint, where $\Tr_B$ is the canonical trace on $B$ (so in
particular we get $\rank B$ when $K=0$). We can then show the
analogues of several of the classical results on pressure as
presented in \cite{W}.

If $(A,\alpha)$ is asymptotically abelian with locality we show the
variational principle (\ref{e1.2}) in Section \ref{3}. Furthermore, if
we assume $ht(\alpha)<\infty$, $H\in\A$ and $\phi$ is a $\beta$-equilibrium
state at~$H$, i.e. $\phi$ is $\alpha$-invariant and $P_\alpha(\beta
H)=h_\phi(\alpha)-\beta\phi(H)$, then we show in Section \ref{4} via
a proof modelled on the corresponding proof based on (\ref{e1.1}) in
\cite{BR2} for spin lattice systems, that $\phi$ is a $\beta$-KMS
state with respect to the one-parameter group defined by the derivation
$$
\delta_H(x)=\sum_{j\in\Z}[\alpha^j(H),x], \ \ x\in\A.
$$
In particular, when $H=0$, so $P_\alpha(0)=ht(\alpha)$, we obtain
\begin{equation} \label{e1.3}
ht(\alpha)=\sup_\phi h_\phi(\alpha),
\end{equation}
and if $h_\phi(\alpha)=ht(\alpha)$ then $\phi$ is a trace.

Equation (\ref{e1.3}) is false in general. Indeed in \cite{NST} there
is exhibited a (non asymptotically abelian) automorphism $\alpha$ of
the CAR-algebra for which the trace $\tau$ is the unique invariant
state, $h_\tau(\alpha)=0$, while $ht(\alpha)\ge{1\over2}\log2$.
Furthermore, our assumption of locality is essential to conclude
that $\phi$ is a trace when $h_\phi(\alpha)=ht(\alpha)<\infty$,
even for asymptotically abelian systems. In Example \ref{5.6} we show
that there is a Bogoliubov automorphism $\alpha$ of the even
CAR-algebra which is asymptotically abelian, $ht(\alpha)=0$, while
there are an infinite number of non-tracial $\alpha$-invariant states.

In Section \ref{5} we consider some other examples and special cases.
First we apply our results to C$^*$-algebras associated with certain
AF-groupoids arising naturally from expansive homeomorphisms of
zero-dimensional compact spaces. We show that if $H$ lies in the
diagonal then there is a one-to-one correspondence between equilibrium
states on the algebra and equilibrium measures in the classical sense.
We also consider unique ergodicity for non-abelian systems. If
$(A,\alpha)$
is asymptotically abelian with locality, unique ergodicity turns out
to be of marginal interest. Indeed, the unique invariant state $\tau$
is a trace, and the image of~$A$ in the GNS-representation of~$\tau$ is
abelian.

\medskip
\noindent
{\bf Acknowledgement.} The authors are indebted to A. Connes for
suggesting to us to study the variational principle and equilibrium
states in the setting of asymptotically abelian C$^*$-algebras.

\bigskip\bigskip

\section{Pressure} \label{2}

In order to define pressure of a $C^*$-dynamical system we follow the
setup of Voiculescu \cite{V} for his definition of topological entropy.

Let $A$ be a nuclear $C^*$-algebra with unit and $\alpha$ an
automorphism. Let $\CPA(A)$ denote the set of triples
$(\rho,\psi,B)$, where $B$ is a finite dimensional $C^{*}$-algebra,
and $\rho\colon A\to B$, $\psi\colon B\to A$ are unital completely
positive maps. ${\cal P}_f(A)$ denotes the family of finite subsets
of~$A$. For $\delta>0$, $\omega\in{\cal P}_f(A)$ and $H\in A_{sa}$, put
$$
P(H,\omega;\delta)=\inf\{\log\Tr_B(e^{-\rho(H)})\,|\,
   (\rho,\psi,B)\in\CPA(A),\,||(\psi\circ\rho)(x)-x||<\delta\ \forall
   x\in\omega\},
$$
where $\Tr_B$ is the canonical trace on $B$, i.e. the trace which
takes the value $1$ on each minimal projection. Then set
\begin{eqnarray*}
P_\alpha(H,\omega;\delta)
   &=&\limsup_{n\to\infty}{1\over n}
         P\left(\sum^{n-1}_{j=0}\alpha^j(H),
            \bigcup^{n-1}_{j=0}\alpha^j(\omega);\delta\right),\\
P_\alpha(H,\omega)&=&\sup_{\delta>0}P_\alpha(H,\omega;\delta).
\end{eqnarray*}

\begin{definition} \label{2.1}
\rm The {\it pressure of~$\alpha$ at $H$} is
$$
P_\alpha(H)=\sup_{\omega\in{\cal P}_f(A)}
   P_\alpha(H,\omega).
$$
\end{definition}

We have chosen the minus sign $e^{-\rho(H)}$ in the definition of
$P(H,\omega;\delta)$ because of its use in physical applications,
see \cite{BR2}, rather than the plus sign used in ergodic theory, see
\cite{W}.

It is easy to see that if $\omega_1\subset\omega_2\subset\ldots$ is
an increasing sequence of finite subsets of~$A$ such
that the linear span of~$\cup_n\omega_n$ is dense in~$A$, then
$P_\alpha(H)=\lim_nP_\alpha(H,\omega_n)$. If $H=0$ the pressure
coincides with the topological entropy of Voiculescu \cite{V}.

Recall that $h_\phi(\alpha)$ denotes the CNT-entropy of~$\alpha$ with
respect to an $\alpha$-invariant state $\phi$ of~$A$.

\begin{proposition} \label{2.2}
For any $\alpha$-invariant state $\phi$ of~$A$ we have
$P_\alpha(H)\ge h_\phi(\alpha)-\phi(H)$.
\end{proposition}

\noindent
{\it Proof.} The proof is a rewording of the proof of
\cite[Proposition~4.6]{V}. Let $N$ be a finite dimensional
C$^*$-algebra, and
$\gamma\colon N\to A$ a unital completely positive map. Let
$\omega\in{\cal P}_f(A)$ be such that $H\in\omega$ and $\gamma(\{x\in
N\,|\,||x||\le1\})$ is contained in the convex hull of~$\omega$. Then
if $(\rho,\psi,B)\in\CPA(A)$ and
$$
||(\psi\circ\rho)(a)-a||<\delta\ \hbox{for}\
a\in\cup^{n-1}_{j=0}\alpha^j(\omega),
$$
we obtain by \cite[Proposition~IV.3]{CNT}
$$
|H_\phi(\{\alpha^j\circ\gamma\}_{0\le j\le n-1})
   -H_\phi(\{\psi\circ\rho\circ\alpha^j\circ\gamma\}_{0\le j\le
   n-1})|<n\eps,
$$
where $\eps=\eps(\delta)\to0$ as $\delta\to0$. If $K\in B_{sa}$ and
$\theta$ is a state on $B$ then $\log\Tr_B(e^{-K})\ge
S(\theta)-\theta(K)$, hence
\begin{eqnarray*}
H_\phi(\{\psi\circ\rho\circ\alpha^j\circ\gamma\}_{0\le j\le n-1})
   &\le&H_\phi(\psi)\le S(\phi\circ\psi)\\
   &\le&\log\Tr_B\left(e^{-\rho\left(
              \sum^{n-1}_{j=0}\alpha^j(H)\right)}\right)
           +(\phi\circ\psi)\left(\rho
      \left(\sum\nolimits^{n-1}_{j=0}\alpha^j(H)\right)\right)\\
   &\le&\log\Tr_B\left(e^{-\rho\left(
              \sum^{n-1}_{j=0}\alpha^j(H)\right)}\right)
           +n\phi(H)+n\delta.
\end{eqnarray*}
Thus
$$
{1\over n}H_\phi(\{\alpha^j\circ\gamma\}_{0\le j\le n-1})
   \le{1\over n}P\left(\sum^{n-1}_{j=0}\alpha^j(H),
         \bigcup^{n-1}_{j=0}\alpha^j(\omega);\delta\right)+\phi(H)
      +\delta+\eps.
$$
It follows that $h_\phi(\gamma;\alpha)\le
P_\alpha(H,\omega)+\phi(H)$, hence
$h_\phi(\alpha)\le P_\alpha(H)+\phi(H)$.
\endverif

\begin{remark} \label{2.3}
\rm If $A$ is abelian, $A=C(X)$, and $P^{cl}_\alpha(H)$ denotes the
pressure
as defined in \cite{W}, then $P_\alpha(H)=P^{cl}_\alpha(-H)$. The
inequality '$\le$' can be proved just the same as in
\cite[Proposition~4.8]{V}. The converse inequality follows from
Proposition~\ref{2.2} and the classical variational principle. In the
AF-case however, i.e. when $X$ is zero-dimensional, it is easy to
give a direct proof. Indeed, if in the proof of Proposition~\ref{2.2}
$N$ was the subalgebra of~$A$ corresponding to a clopen partition $P$
of~$X$, then we could conclude that
\begin{equation} \label{e2.1}
{1\over n}H_\phi(\{\alpha^j(N)\}_{0\le j\le n-1})
   \le{1\over n}P\left(\sum^{n-1}_{j=0}\alpha^j(H),
         \bigcup^{n-1}_{j=0}\alpha^j(\omega);\delta\right)
       +{1\over n}\phi\left(\sum^{n-1}_{j=0}\alpha^j(H)\right)
       +\delta+\eps
\end{equation}
for any (not necessarily $\alpha$-invariant) state $\phi$ of~$A$,
with $\eps$ independent of~$\phi$. Let $T$ be the homeomorphism
corresponding to $\alpha$, so that $\alpha(f)=f\circ T$. Suppose
the points $x_1,\ldots, x_m$ lie in different elements of the partition
$\vee^{n-1}_{j=0}T^{-j}P$. Then inequality (\ref{e2.1}) for the
measure
$$
\phi=\left(\sum_ie^{-(S_nH)(x_i)}\right)^{-1}
      \sum_ie^{-(S_nH)(x_i)}\delta_{x_i},
$$
where $(S_nH)(x)=\sum^{n-1}_{j=0}H(T^jx)$, means that
$$
{1\over n}\log\sum_ie^{-(S_nH)(x_i)}
   \le{1\over n}P\left(\sum^{n-1}_{j=0}\alpha^j(H),
         \bigcup^{n-1}_{j=0}\alpha^j(\omega);\delta\right)
       +\delta+\eps.
$$
Recalling the definition of pressure \cite[Definition~9.7]{W}, we see
that $P^{cl}_\alpha(-H)\le P_\alpha(H)$.
\end{remark}

We list some properties of the function $H\mapsto P_\alpha(H)$
on~$A_{sa}$.

\begin{proposition} \label{2.4}
The following properties are satisfied by $P_\alpha$ for $H,K\in A_{sa}$.
\begin{list}{}{}
\item{\rm(i)} If $H\le K$ then $P_\alpha(H)\ge P_\alpha(K)$.
\item{\rm(ii)} $P_\alpha(H+c1)=P_\alpha(H)-c$, $c\in\R$.
\item{\rm(iii)} $P_\alpha(H)$ is either infinite for all $H$ or is finite
valued.
\item{\rm(iv)} If $P_\alpha$ is finite valued then
$|P_\alpha(H)-P_\alpha(K)|\le||H-K||$.
\item{\rm(v)} For $k\in\N$,
$P_{\alpha^k}(\sum^{k-1}_{j=0}\alpha^j(H))=kP_\alpha(H)$.
\item{\rm(vi)} $P_\alpha(H+\alpha(K)-K)=P_\alpha(H)$.
\end{list}
\end{proposition}

\noindent
{\it Proof.} (i) Given $H\le K$ take $\omega\in{\cal P}_f(A)$.
If $(\rho,\psi,B)\in\CPA(A)$ we have
$$
\log\Tr_B\left(e^{-\rho\left(\sum^{n-1}_{j=0}\alpha^j(H)\right)}\right)
   \ge\log\Tr_B\left(e^{-\rho\left(
            \sum^{n-1}_{j=0}\alpha^j(K)\right)}\right),
$$
see e.g. \cite[Corollary~3.15]{OP}. Thus (i) follows.

(ii) As in (i) we have
$$
\log\Tr_B\left(e^{-\rho\left(\sum^{n-1}_{j=0}\alpha^j(H+c1)\right)}\right)
   =\log\Tr_B\left(e^{-\rho\left(
         \sum^{n-1}_{j=0}\alpha^j(H)\right)}\right)-nc,
$$
and (ii) follows.

(iii) By (i) and (ii) we have
$$
P_\alpha(H)\ge P_\alpha(||H||)=P_\alpha(0)-||H||=ht(\alpha)-||H||,
$$
and similarly $P_\alpha(H)\le ht(\alpha)+||H||$. Thus (iii) follows.

(iv) For any $(\rho,\psi,B)\in\CPA(A)$ we have by the
Peierls-Bogoliubov inequality \cite[Corollary~3.15]{OP}
$$
\left|{1\over n}\log\Tr_B\left(e^{-\rho\left(
      \sum^{n-1}_{j=0}\alpha^j(H)\right)}\right)
   -{1\over n}\log\Tr_B\left(e^{-\rho\left(
           \sum^{n-1}_{j=0}\alpha^j(K)\right)}\right)\right|\le||H-K||.
$$
Thus
$$
P_\alpha(H,\omega;\delta)-||H-K||\le P_\alpha(K,\omega;\delta)
   \le P_\alpha(H,\omega;\delta)+||H-K||
$$
for any $\omega\in{\cal P}_f(A)$. Thus (iv)
follows.

(v) Let $(\rho,\psi,B)\in\CPA(A)$ and
$\omega\in{\cal P}_f(A)$. Given $n\in\N$ choose $m\in\N$ such that $mk\le
n<(m+1)k$. Set $H_k=\sum^{k-1}_{j=0}\alpha^j(H)$ and
$\omega_k=\cup^{k-1}_{j=0}\alpha^j(\omega)$. Then
\begin{eqnarray*}
\log\Tr_B\left(e^{-\rho\left(\sum^{n-1}_{j=0}\alpha^j(H)\right)}\right)
  &\ge&\log\Tr_B\left(e^{-\rho\left(
     \sum^{mk-1}_{j=0}\alpha^j(H)\right)-k||H||}\right)\\
  & = &\log\Tr_B\left(e^{-\rho\left(\sum^{m-1}_{j=0}\alpha^{jk}(H_k)\right)}
          \right)-k||H||.
\end{eqnarray*}
Similarly
$$
\log\Tr_B\left(e^{-\rho\left(\sum^{n-1}_{j=0}\alpha^j(H)\right)}\right)
 \le\log\Tr_B\left(e^{-\rho\left(\sum^m_{j=0}\alpha^{jk}(H_k)\right)}
       \right)+k||H||.
$$
Since $\cup^{m-1}_{j=0}\alpha^{jk}(\omega_k)\subset
\cup^{n-1}_{j=0}\alpha^j(\omega)\subset
\cup^{m}_{j=0}\alpha^{jk}(\omega_k)$, it follows that
$$
P_\alpha(H,\omega;\delta)={1\over k}P_{\alpha^k}
   \left(\sum^{k-1}_{j=0}\alpha^j(H),
      \bigcup^{k-1}_{j=0}\alpha^j(\omega);\delta\right),
$$
and hence $P_\alpha(H)={1\over k}P_{\alpha^k}(\sum^{k-1}_{j=0}\alpha^j(H))$.

(vi) Set $H_k=\sum^{k-1}_{j=0}\alpha^j(H)$ and
$H'_k=\sum^{k-1}_{j=0}\alpha^j(H+\alpha(K)-K)=H_k+\alpha^k(K)-K$. Then
by (iv) and (v) we have
$$
|P_\alpha(H)-P_\alpha(H+\alpha(K)-K)|
   ={1\over k}|P_{\alpha^k}(H_k)-P_{\alpha^k}(H'_k)|
   \le{2||K||\over k}.
$$
Thus (vi) follows.
\endverif

The next result is the analogue of \cite[Theorem~9.12]{W}, see also
\cite{R}.

\begin{proposition} \label{2.5}
Suppose $ht(\alpha)<\infty$. Let $\phi$ be a self-adjoint linear
functional on $A$. Then $\phi$ is an $\alpha$-invariant state if and
only if $-\phi(H)\le P_\alpha(H)$ for all $H\in A_{sa}$.
\end{proposition}

\noindent
{\it Proof.} If $\phi$ is an $\alpha$-invariant state then by
Proposition~\ref{2.2}
$$
-\phi(H)\le P_\alpha(H)-h_\phi(\alpha)\le P_\alpha(H)
$$
for all $H\in A_{sa}$.

Conversely if $-\phi(H)\le P_\alpha(H)$ for all $H\in A_{sa}$ then by
Proposition~\ref{2.4}(vi)
$$
-\phi(\alpha(H)-H)=-{1\over n}\phi(\alpha(nH)-nH)
   \le{1\over n}P_\alpha(\alpha(nH)-nH)
   ={1\over n}P_\alpha(0)
   ={1\over n}ht(\alpha)\mathop{\longrightarrow}_{n\to\infty}0.
$$
Applying this also to $-H$ we see that $\phi$ is $\alpha$-invariant.
Furthermore, by Proposition~\ref{2.4}(i),(ii)
$$
-\phi(H)=-{1\over n}\phi(nH)
   \le{1\over n}P_\alpha(nH)
   \le{1\over n}ht(\alpha)+||H||\mathop{\longrightarrow}_{n\to\infty}||H||,
$$
so that $||\phi||\le1$. For $c\in\R$ we have
$$
-c\phi(1)\le P_\alpha(c1)=ht(\alpha)-c.
$$
Hence $\phi(1)=1$, and $\phi$ is a state.
\endverif

\begin{definition} \label{2.6}
\rm We say $\phi$ is an {\it equilibrium state at $H$} if
$$
P_\alpha(H)=h_\phi(\alpha)-\phi(H),
$$
hence by Proposition~\ref{2.2}
$$
h_\phi(\alpha)-\phi(H)=\sup_\psi(h_\psi(\alpha)-\psi(H)),
$$
where the sup is taken over all $\alpha$-invariant states.
\end{definition}

Recall that if $F$ is a real convex continuous function on a real Banach
space $X$, then a linear functional $f$ on $X$ is called a tangent
functional to the graph of~$F$ at the point $x_0\in X$ if
$$
F(x_0+x)-F(x)\ge f(x)\ \forall x\in X.
$$
In the sequel we will identify self-adjoint linear functionals on $A$
with real linear functionals on~$A_{sa}$.

The next proposition is the analogue of \cite[Theorem~9.14]{W}.

\begin{proposition} \label{2.7}
Suppose $ht(\alpha)<\infty$ and the pressure is a convex function
on~$A_{sa}$. Then
\begin{list}{}{}
\item{\rm (i)} if $\phi$ is an equilibrium state at $H$ then $-\phi$
is a tangent functional for the pressure at $H$;
\item{\rm (ii)} if $-\phi$ is a tangent functional for
the pressure at $H$ then $\phi$ is an $\alpha$-invariant state.
\end{list}
\end{proposition}

\noindent
{\it Proof.} (i) Let $K\in A_{sa}$. Then by Proposition~\ref{2.2}
$$
P_\alpha(H+K)-P_\alpha(H)\ge(h_\phi(\alpha)-\phi(H+K))
   -(h_\phi(\alpha)-\phi(H))=-\phi(K),
$$
so $-\phi$ is a tangent functional.

(ii) If $K\in A_{sa}$ then by Proposition~\ref{2.4}(vi)
$$
-\phi(\alpha(K)-K)\le P_\alpha(H+\alpha(K)-K)-P_\alpha(H)=0.
$$
Applying this also to $-K$ we see that $\phi$ is $\alpha$-invariant.
Now note that $||\phi||\le1$ by Proposition~\ref{2.4}(iv). By
Proposition~\ref{2.4}(ii) we have also $-c\ge-c\phi(1)$ for any
$c\in\R$. Hence $\phi(1)=1$, and $\phi$ is a state.
\endverif

\bigskip\bigskip

\section{The variational principle} \label{3}

We shall prove the variational principle for the following class of
C$^*$-dynamical systems.

\begin{definition} \label{3.1}
\rm A unital C$^*$-dynamical system $(A,\alpha)$ is called
{\it asymptotically abelian with locality} if there is a dense
$\alpha$-invariant $*$-subalgebra $\cal A$ of~$A$ such that for each
pair $a,b\in\cal A$ the C$^*$-algebra generated by $a$ and $b$ is
finite dimensional, and for some $p=p(a,b)\in\N$ we have $[\alpha^j(a),b]=0$
whenever $|j|\ge p$.
\end{definition}

We call elements of~$\cal A$ for {\it local operators} and finite
dimensional C$^*$-subalgebras of~$\cal A$ for {\it local algebras}.

Note that since we may add the identity operator to $\cal A$, we may
assume that $1\in\cal A$. Since each finite dimensional
C$^*$-algebra is singly generated, an easy induction argument shows
that the C$^*$-algebra generated by a finite set of local operators
is finite dimensional. In particular,
$A$ is an AF-algebra. Note also that another
easy induction argument shows that for each local algebra $N$ there
is $p\in\N$ such that $[\alpha^j(a),b]=0$ for all $a,b\in N$
whenever $|j|\ge p$.

\begin{theorem} \label{3.2}
Let $(A,\alpha)$ be a unital separable C$^*$-dynamical system which
is asymptotically abelian with locality. Let $H\in A_{sa}$. Then
$$
P_\alpha(H)=\sup_\phi(h_\phi(\alpha)-\phi(H)),
$$
where the sup is taken over all $\alpha$-invariant states of~$A$. In
particular, the topological entropy
$$
ht(\alpha)=\sup_\phi h_\phi(\alpha).
$$
\end{theorem}

Consider first the case when there exists a finite dimensional
C$^*$-subalgebra $N$ of~$A$ such that $H\in N$, $\alpha^j(N)$
commutes with $N$ for $j\ne0$, $\vee_{j\in\Z}\alpha^j(N)=A$.

\begin{lemma} \label{3.3}
Under the above assumptions there exists an $\alpha$-invariant
state $\phi$ such that
$$
P_\alpha(H)=h_\phi(\alpha)-\phi(H)=\lim_{n\to\infty}{1\over n}
   \log\Tr_{\vee^{n-1}_{j=0}\alpha^j(N)}
      \left(e^{-\sum^{n-1}_{j=0}\alpha^j(H)}\right).
$$
\end{lemma}

\noindent
{\it Proof.} First note that if $A_1$ and $A_2$ are commuting finite
dimensional C$^*$-algebras, and $a_i\in A_i$, $a_i\ge0$, $i=1,2$,
then
$$
\Tr_{A_1\vee A_2}(a_1a_2)\le\Tr_{A_1}(a_1)\Tr_{A_2}(a_2),
$$
since if $p_i$ is a minimal projection in~$A_i$, $i=1,2$, then
$p_1p_2$ is either zero or minimal in~$A_1\vee A_2$. Hence the limit
in the formulation of the lemma really exists. We denote it by $\tilde
P_\alpha(H)$. It is easy to see that $\tilde P_\alpha(H)\ge P_\alpha(H)$.

For each $n\in\Z$ let $N_n$ be a copy of~$N$. Consider the
C$^*$-algebra $M=\otimes_{n\in\Z}N_n$. Let $\beta$ be the shift to
the right on $M$, and $\pi\colon M\to A$ the homomorphism which
intertwines $\beta$ with $\alpha$, and identifies $N_0$ with $N$.
Set $I=\Ker\,\pi$. For each $n\in\N$ let
$$
M_n=N_0\otimes\ldots\otimes N_{n-1},\ \
I_n=I\cap M_n, \ \ \pi_n=\pi|_{M_n}.
$$
Identifying $M$ with $M^{\otimes\Z}_n$ consider the
$\beta^n$-invariant state $\psi_n=\otimes(f_n\circ\pi_n)$ on $M$,
where $f_n$ is the state on $\vee^{n-1}_{j=0}\alpha^j(N)$ with
density operator
$$
\left(\Tr_{\vee^{n-1}_{j=0}\alpha^j(N)}e^{-\sum^{n-1}_{j=0}\alpha^j(H)}
      \right)^{-1}e^{-\sum^{n-1}_{j=0}\alpha^j(H)}.
$$
Set $\displaystyle\phi_n={1\over n}\sum^{n-1}_{j=0}\psi_n\circ\beta^j$.
Then $\phi_n$ is $\beta$-invariant.
Using concavity of entropy we obtain
\begin{eqnarray*}
h_{\phi_n}(\beta)
 &=&{1\over n}h_{\phi_n}(\beta^n)
       \ge{1\over n^2}\sum^{n-1}_{j=0}h_{\psi_n\circ\beta^j}(\beta^n)
       ={1\over n}h_{\psi_n}(\beta^n)
       ={1\over n}S(f_n\circ\pi_n)
       ={1\over n}S(f_n)\\
 &=&{1\over n}\log\Tr_{\vee^{n-1}_{j=0}\alpha^j(N)}\left(
          e^{-\sum^{n-1}_{j=0}\alpha^j(H)}\right)
       +{1\over n}f_n\left(\sum^{n-1}_{j=0}\alpha^j(H)\right)\\
 &\ge&\tilde P_\alpha(H)
       +{1\over n}f_n\left(\sum^{n-1}_{j=0}\alpha^j(H)\right)\\
 &=&\tilde P_\alpha(H)
       +{1\over n}\psi_n\left(\sum^{n-1}_{j=0}\beta^j(H)\right)
       =\tilde P_\alpha(H)+\phi_n(H).
\end{eqnarray*}
Let $\tilde\phi$ be any weak$^{*}$ limit point of the sequence
$\{\phi_n\}_n$. Then $\tilde\phi$ is $\beta$-invariant. Let $B$ be a masa
in~$N_0$ containing $H$. Then $B$ is
in the centralizer of the state $\phi_n$, hence
$$
h_{\phi_n}(\beta)=h_{\phi_n}(B;\beta).
$$
Since the mapping $\psi\mapsto h_\psi(B;\beta)$ is upper
semicontinuous, we conclude that
$$
h_{\tilde\phi}(\beta)\ge\tilde P_\alpha(H)+\tilde\phi(H).
$$
Now note that $\tilde\phi$ is zero on $I$. Indeed, if $x\in I_n$
then $\beta^j(x)\in I_m$ for $j=0,\ldots,m-n$ and $m\ge n$, whence
$$
|\phi_m(x)|\le{1\over m}\sum^{m-1}_{j=m-n+1}|(\psi_m\circ\beta^j)(x)|
   \le{n-1\over m}||x||,
$$
so $\tilde\phi(x)=0$. Thus $\tilde\phi$ defines a state $\phi$ on
$A$. We have
$$
h_\phi(\alpha)=h_{\tilde\phi}(\beta)\ge\tilde P_\alpha(H)+\phi(H),
$$
where the first equality follows from \cite[Theorem~VII.2]{CNT}.
Since by Proposition~\ref{2.2},
$h_\phi(\alpha)-\phi(H)\le P_\alpha(H)\le\tilde
P_\alpha(H)$, the proof of the lemma is complete.
\endverif

We shall reduce the general case to the case considered above by
replacing $\alpha$ by its powers. For this suppose that $N$ is a
local subalgebra of~$A$, and $H\in N$. Choose $p$ such that
$\alpha^j(N)$ commutes with $N$ whenever $|j|\ge p$. For $k\ge p$ set
$M_k=\vee^{k-p}_{j=0}\alpha^j(N)$, $H_k=\sum^{k-p}_{j=0}\alpha^j(H)$.
Then $H_k\in M_k$, and $\alpha^{jk}(M_k)$ commutes with $M_k$ for $j\ne0$.

\begin{lemma} \label{3.4}
For any finite subset $\omega$ of~$N$ we have
$$
P_\alpha(H,\omega)\le\liminf_{k\to\infty}{1\over k}
   P_{\alpha^k|\vee_{j\in\Z}\alpha^{jk}(M_k)}(H_k).
$$
\end{lemma}

\noindent
{\it Proof.} The idea of the proof is to reduce to the situation of
Lemma~\ref{3.3} by showing that the contribution of the indicies in the
intervals $[jk-p+1,jk-1]$, $j\in\N$, becomes negligible for large $k$.

 Fix $\delta>0$. Choose $m_0\in\N$ such that
$$
{2(p-1)||a||\over m_0}<\delta\ \ \hbox{for}\ \ a\in\omega.
$$
Take any $k\ge m_0+p$. Let $n\in\N$. Then $(m-1)k\le n<mk$ for
some $m\in\N$. Set $B_0=\vee^m_{j=0}\alpha^{jk}(M_k)$ and
$$
B=\underbrace{B_0\oplus\ldots\oplus B_0}_{m_0}.
$$
Choose a conditional expectation $E\colon A\to B_0$, and define unital
completely positive mappings $\psi\colon B\to A$ and $\rho\colon A\to
B$ as follows:
$$
\psi(b_1,\ldots,b_{m_0})={1\over m_0}
   \sum^{m_0}_{i=1}\alpha^{-i+1}(b_i),
$$
$$
\rho(a)=(E(a),(E\circ\alpha)(a),\ldots,(E\circ\alpha^{m_0-1})(a)).
$$
For any $a\in A$ we have
$$
||(\psi\circ\rho)(a)-a||\le{2||a||\over m_0}
   \#\{0\le i\le m_0-1\,|\,\alpha^i(a)\notin B_0\},
$$
where $\#S$ means the cardinality of a set $S$.
Let $a=\alpha^l(b)$ for some $b\in\omega$ and $l$, $0\le l\le n-1$. Then
$l=jk+r$ for some $j$ and $r$, $0\le j\le m-1$, $0\le r<k$. Since
$m_0\le k-p$, the interval $[l,l+m_0-1]$ is contained in
$[jk,(j+1)k+k-p]$. But for
$i\in[jk,(j+1)k+k-p]\backslash[jk+k-p+1,(j+1)k-1]$ we have
$\alpha^i(N)\subset B_0$, so
$\#\{0\le i\le m_0-1\,|\,\alpha^i(a)\notin B_0\}\le p-1$, and
$$
||(\psi\circ\rho)(a)-a||\le{2(p-1)||b||\over m_0}<\delta.
$$
Hence
$$
P\left(\sum^{n-1}_{j=0}\alpha^j(H),\bigcup^{n-1}_{j=0}\alpha^j(\omega);
      \delta\right)
   \le\log\Tr_B\left(e^{-\rho\left(
            \sum^{n-1}_{j=0}\alpha^j(H)\right)}\right).
$$

Now note that for $0\le i\le m_0-1$ the sets $X_i=[i,i+n-1]$ and
$X=\cup^m_{j=0}[jk,jk+k-p]$ are contained in~$Y=[0,mk+k-p]$, so
\begin{eqnarray*}
\#(X_i\bigtriangleup X)
 &\le&\#(Y\backslash X_i)+\#(Y\backslash X)=(mk+k-p+1-n)+m(p-1)\\
 &\le&mk+k-p+1-(m-1)k+m(p-1)\le mp+2k.
\end{eqnarray*}
When $j\in X_i\cap X$, $\alpha^j(H)\in B_0$. Hence
$$
\left\|(E\circ\alpha^i)\left(\sum^{n-1}_{j=0}\alpha^j(H)\right)
   -\sum^m_{j=0}\alpha^{jk}(H_k)\right\|
   =\left\|E\left(\sum_{j\in X_i}\alpha^j(H)\right)
      -\sum_{j\in X}\alpha^j(H)\right\|
   \le(mp+2k)||H||.
$$
By the Peierls-Bogoliubov inequality we obtain
$$
\Tr_{B_0}\left(e^{-(E\circ\alpha^i)\left(\sum^{n-1}_{j=0}\alpha^j(H)\right)}
      \right)
   \le e^{(mp+2k)||H||}
         \Tr_{B_0}\left(e^{-\sum^m_{j=0}\alpha^{jk}(H_k)}\right),
$$
so
$$
\Tr_B\left(e^{-\rho\left(\sum^{n-1}_{j=0}\alpha^j(H)\right)}\right)
   \le m_0e^{(mp+2k)||H||}
         \Tr_{B_0}\left(e^{-\sum^m_{j=0}\alpha^{jk}(H_k)}\right).
$$
Taking the $\log$, dividing by $n$, and letting $n\to\infty$, we obtain
\begin{eqnarray*}
P_\alpha(H,\omega;\delta)
 &\le&{p||H||\over k}+{1\over k}\lim_{m\to\infty}{1\over m}\log
         \Tr_{\vee^{m-1}_{j=0}\alpha^{jk}(M_k)}\left(
    e^{-\sum^{m-1}_{j=0}\alpha^{jk}(H_k)}\right)\\
 &=&{p||H||\over k}+{1\over k}
       P_{\alpha^k|\vee_{j\in\Z}\alpha^{jk}(M_k)}(H_k),
\end{eqnarray*}
where the last equality follows from Lemma~\ref{3.3}.
\endverif

We shall need also the following
\begin{lemma} \label{3.5}
Let $(A,\alpha)$ be a C$^*$-dynamical system with A nuclear, $B$
an $\alpha$-invariant C$^*$-subalgebra of~$A$, $\phi$ an $\alpha$-invariant
state on $B$. Then for any $\eps>0$ there exists an $\alpha$-invariant
state $\psi$ on $A$ such that
$$
\psi|_B=\phi\ \ \hbox{and}\ \ h_\psi(\alpha)>h_\phi(\alpha|_B)-\eps.
$$
\end{lemma}

\noindent
{\it Proof.} Since the Sauvageot-Thouvenot entropy is not less than
the CNT-entropy for general C$^*$-systems, there exist a
commutative C$^*$-dynamical system $(C,\beta,\mu)$, an
$(\alpha\otimes\beta)$-invariant state $\lambda$ on $B\otimes C$, and
a finite dimensional subalgebra $P$ of~$C$ such that
$\lambda|_B=\phi$, $\lambda|_C=\mu$ and
$$
h_\phi(\alpha|_B)<H_\mu(P,P^-)-H_\lambda(P|B)+\eps,
$$
see \cite{ST} for notations. Extend $\lambda$ to an
$(\alpha\otimes\beta)$-invariant state $\Lambda$ on $A\otimes C$, and
set $\psi=\Lambda|_A$. Since the conditional entropy $H_\lambda(P|B)$
is decreasing in second variable, and ST-entropy coincides with
CNT-entropy for nuclear algebras, we have
$$
h_\psi(\alpha)\ge H_\mu(P,P^-)-H_\Lambda(P|A)
   \ge H_\mu(P,P^-)-H_\lambda(P|B)>h_\phi(\alpha|_B)-\eps.
$$
\endveriff

\noindent
{\it Proof of Theorem~\ref{3.2}.} The inequality '$\ge$' has been
proved in Proposition~\ref{2.2}. Since the pressure is continuous by
Proposition~\ref{2.4}(iv), to prove the converse inequality it suffices
to consider local~$H$. Then
by Lemma~\ref{3.4} we have only to show that if $H$ is contained in a
local algebra $N$ then
$$
\sup_\phi(h_\phi(\alpha)-\phi(H))\ge\liminf_{k\to\infty}{1\over k}
   P_{\alpha^k|\vee_{j\in\Z}\alpha^{jk}(M_k)}(H_k).
$$
By Lemma~\ref{3.3}, for each $k\in\N$ there exists an
$\alpha^k$-invariant state $\psi_k$ on
$\vee_{j\in\Z}\alpha^{jk}(M_k)$ such that
$$
h_{\psi_k}(\alpha^k|_{\vee_{j\in\Z}\alpha^{jk}(M_k)})-\psi_k(H_k)
   =P_{\alpha^k|\vee_{j\in\Z}\alpha^{jk}(M_k)}(H_k).
$$
By Lemma~\ref{3.5} we may extend $\psi_k$ to an $\alpha^k$-invariant
state $\tilde\phi_k$ on $A$ such that
$$
h_{\tilde\phi_k}(\alpha^k)
   \ge h_{\psi_k}(\alpha^k|_{\vee_{j\in\Z}\alpha^{jk}(M_k)})-1.
$$
Set $\displaystyle\phi_k={1\over
k}\sum^{k-1}_{j=0}\tilde\phi_k\circ\alpha^j$. Then as in the proof of
Lemma~\ref{3.3}
$$
h_{\phi_k}(\alpha)
   \ge{1\over k}h_{\tilde\phi_k}(\alpha^k)
   \ge{1\over k}h_{\psi_k}(\alpha^k|_{\vee_{j\in\Z}\alpha^{jk}(M_k)})
         -{1\over k}.
$$
Since
$$
\phi_k(H)={1\over k}\sum^{k-1}_{j=0}\tilde\phi_k(\alpha^j(H))
   \le{1\over k}\tilde\phi_k(H_k)+{p-1\over k}||H||
   ={1\over k}\psi_k(H_k)+{p-1\over k}||H||,
$$
we get
\begin{eqnarray*}
h_{\phi_k}(\alpha)-\phi_k(H)
 &\ge&{1\over k}h_{\psi_k}(\alpha^k|_{\vee_{j\in\Z}\alpha^{jk}(M_k)})
         -{1\over k}\psi_k(H_k)-{1+(p-1)||H||\over k}\\
 &=&{1\over k}P_{\alpha^k|\vee_{j\in\Z}\alpha^{jk}(M_k)}(H_k)
       -{1+(p-1)||H||\over k},
\end{eqnarray*}
and the proof is complete.
\endverif

\begin{corollary} \label{3.6}
With our assumptions the pressure is a convex function of~$H$.
\end{corollary}

\noindent
{\it Proof.} Use the affinity of the function
$H\mapsto h_\phi(\alpha)-\phi(H)$.
\endverif

\begin{corollary} \label{3.7}
If $(A_1,\alpha_1)$ and $(A_2,\alpha_2)$ are asymptotically abelian
systems with locality then
$$
ht(\alpha_1\otimes\alpha_2)=ht(\alpha_1)+ht(\alpha_2).
$$
\end{corollary}

\noindent
{\it Proof.} If $\phi_i$ is an $\alpha_i$-invariant state, $i=1,2$,
then by \cite[Lemma~3.4]{SV} and \cite[Propositions~4.6 and~4.9]{V},
$$
h_{\phi_1}(\alpha_1)+h_{\phi_2}(\alpha_2)
   \le h_{\phi_1\otimes\phi_2}(\alpha_1\otimes\alpha_2)
   \le ht(\alpha_1\otimes\alpha_2)
   \le ht(\alpha_1)+ht(\alpha_2).
$$
Taking the sup over $\phi_i$ we get the conclusion.
\endverif

\bigskip\bigskip

\section{KMS-states} \label{4}

By Corollary \ref{3.6} and Proposition~\ref{2.7} it follows that if
$(A,\alpha)$ is asymptotically abelian with locality and
$ht(\alpha)<\infty$, then for every equilibrium state $\phi$ at $H$,
$-\phi$ is a tangent functional for the pressure $P_\alpha$ at $H$.
Furthermore, if $\omega$ is a tangent functional for $P_\alpha$ at $H$
then $-\omega$ is an $\alpha$-invariant state.

If $H$ is local and $I\subset\Z$ is a subset then the derivation
$$
\delta_{H,I}(x)=\sum_{j\in I}[\alpha^j(H),x],\ \ x\in{\cal A},
$$
defines a strongly continuous one-parameter automorphism group
$\sigma^{H,I}_t=\exp(it\delta_{H,I})$ of~$A$ (see \cite[Theorem~6.2.6
and~Example 6.2.8]{BR2}). We shall mainly be
concerned with the case $I=\Z$, and will write
$\delta_H=\delta_{H,\Z}$, $\sigma^H_t=\sigma^{H,\Z}_t$. Recall that a
state $\phi$ is a $(\sigma^H_t,\beta)$-KMS state if
$\phi(ab)=\phi(b\sigma^H_{i\beta}(a))$ for $\sigma^H_t$-analytic
elements $a,b\in A$.

We say that an $\alpha$-invariant state $\phi$ is an {\it equilibrium state
at $H$ at inverse temperature $\beta$} if
$$
P_\alpha(\beta H)=h_\phi(\alpha)-\beta\phi(H).
$$
By Theorem~\ref{3.2}, for systems which are asymptotically abelian with
locality, this is equivalent to
$$
h_\phi(\alpha)-\beta\phi(H)=\sup_\psi(h_\psi(\alpha)-\beta\psi(H)).
$$
The main result in this section is

\begin{theorem} \label{4.1}
Suppose a unital separable C$^*$-dynamical system $(A,\alpha)$ is
asymptotically abelian with locality,
and $ht(\alpha)<\infty$. If $H$ is a local self-adjoint operator in~$A$
and
$\phi$ is an equilibrium state at $H$ at inverse temperature~$\beta$,
then $\phi$ is a $(\sigma^H_t,\beta)$-KMS state.
In particular, if $ht(\alpha)=h_\phi(\alpha)$ then $\phi$ is a trace.
\end{theorem}

In order to prove the theorem we may replace $H$ by $\beta H$ and show
that $\phi$ is a $(\sigma^H_t,1)$-KMS state. We shall prove the
following more general result.

\begin{theorem} \label{4.2}
If $-\phi$ is a tangent functional for $P_\alpha$ at $H$ then $\phi$
is a $(\sigma^H_t,1)$-KMS state.
\end{theorem}

We shall need an explicit formula for the pressure, which is a
consequence of our proof of the variational principle.

\begin{lemma} \label{4.3}
Let $N$ be a local algebra. Then there exist a sequence $\{A_n\}_n$
of local algebras containing $N$ and three sequences $\{p_n\}_n$,
$\{m_n\}_n$, $\{k_n\}_n$ of positive integers such that
\begin{list}{}{}
\item{\rm(i)} $\alpha^p(A_n)$ commutes with $A_n$ whenever $|p|\ge p_n$;
\item{\rm(ii)} $\displaystyle{p_n\over k_n}\to0$ as $n\to\infty$;
\item{\rm(iii)}
$\displaystyle P_\alpha(H)=\lim_{n\to\infty}{1\over k_nm_n}
   \log\Tr_{\vee_{j\in I_n}\alpha^j(A_n)}
      \left(e^{-\sum_{j\in I_n}\alpha^j(H)}\right)$ for all $H\in
N_{sa}$, where \newline
$\displaystyle I_n=\bigcup^{m_n-1}_{j=0}[jk_n,jk_n+k_n-p_n]$.
\end{list}
\end{lemma}

\noindent
{\it Proof.} Let $\{A_n\}_n$ be an increasing sequence of local
algebras containing $N$ such that $\cup_nA_n$ is dense in~$A$,
$\omega_n$ a finite subset of~$A_n$ such that $\span(\omega_n)=A_n$.
Let $\{p_n\}_n$ be a sequence satisfying condition (i). By Lemma~\ref{3.4}
$$
P_\alpha(H,\omega_n)\le\liminf_{k\to\infty}{1\over k}
   P_{\alpha^k|\vee_{j\in\Z}\alpha^{jk}(A_{n,k})}\left(
      \sum^{k-p_n}_{j=0}\alpha^j(H)\right)\ \ \forall H\in N_{sa},
$$
where $A_{n,k}=\vee^{k-p_n}_{j=0}\alpha^j(A_n)$. On the other hand,
by the proof of Theorem~\ref{3.2}
$$
P_\alpha(H)\ge\limsup_{k\to\infty}{1\over k}
   P_{\alpha^k|\vee_{j\in\Z}\alpha^{jk}(A_{n,k})}\left(
      \sum^{k-p_n}_{j=0}\alpha^j(H)\right)\ \ \forall H\in N_{sa}.
$$
Choose a countable dense subset $X$ of~$N_{sa}$. Since
$P_\alpha(H,\omega_n)\nearrow P_\alpha(H)$ for any $H\in X$, we can
find a sequence $\{k_n\}_n$ such that condition (ii) is satisfied
and
$$
P_\alpha(H)=\lim_{n\to\infty}{1\over k_n}
   P_{\alpha^{k_n}|\vee_{j\in\Z}\alpha^{jk_n}(A_{n,k_n})}\left(
      \sum^{k_n-p_n}_{j=0}\alpha^j(H)\right)\ \ \forall H\in X.
$$
Since by Lemma~\ref{3.3}
$$
P_{\alpha^{k_n}|\vee_{j\in\Z}\alpha^{jk_n}(A_{n,k_n})}\left(
      \sum^{k_n-p_n}_{j=0}\alpha^j(H)\right)
   =\lim_{m\to\infty}{1\over m}\log\Tr_{\vee_{j\in I_{n,m}}\alpha^j(A_n)}
      \left(e^{-\sum_{j\in I_{n,m}}\alpha^j(H)}\right),
$$
where $I_{n,m}=\cup^{m-1}_{j=0}[jk_n,jk_n+k_n-p_n]$,
we can choose a sequence $\{m_n\}_n$ such that condition~(iii) is
satisfied for all $H\in X$. But then it is satisfied for all $H\in
N_{sa}$ by Proposition~\ref{2.4}(iv) and the Peierls-Bogoliubov
inequality.
\endverif

Every local operator is analytic for the dynamics, and $\sigma^H_t$
depends continuously on $H$ in a fixed local algebra. More precisely,
we have

\begin{lemma} \label{4.4}
\mbox{\ }

{\rm(i)} The series $\displaystyle\sigma^{H,I}_\beta(a)=\sum^\infty_{n=0}
{(i\beta)^n\over n!}\delta^n_{H,I}(a)$ converges absolutely in norm for any
$\beta\in\C$ and any local operator $a$.

{\rm(ii)} Given a local algebra $N$, $R>0$, $C>0$ and $\eps>0$
there exist $q\in\N$ and $\delta>0$ such that
$$
||\sigma^{H_1,I_1}_\beta(a)-\sigma^{H_2,I_2}_\beta(a)||\le\eps||a||
$$
$\forall a\in N$, $\forall H_1,H_2\in N_{sa}$ with $||H_1||,||H_2||\le
C$ and $||H_1-H_2||<\delta$, $\forall\beta\in\C$ with $|\beta|\le R$,
$\forall I_1,I_2\subset\Z$ with $[-q,q]\subset I_1\cap I_2$.
\end{lemma}

\noindent
{\it Proof.} We shall use the arguments of Araki \cite[Theorem~4.2]{A}.
Let $H$ and $a$ lie in a local algebra~$N$. Choose $p\in\N$ such
that $\alpha^j(N)$ commutes with $N$ for $|j|\ge p$. Then
$$
\delta^m_{H,I}(a)=\sum_{j_1,\ldots,j_m}[\alpha^{j_m}(H),[\ldots,
   [\alpha^{j_1}(H),a]\ldots]],
$$
where the sum is over all $j_1,\ldots,j_m$ such that
\begin{equation} \label{e4.1}
j_k\in[-p,p]\bigcup\left(\bigcup_{l<k}[j_l-p,j_l+p]\right)
\end{equation}
for each $k=1,\ldots,m$. But as was already noted in \cite{GN}
condition (\ref{e4.1}) is equivalent to
$$
[j_k,j_k+p]\cap\left([0,p]\bigcup\left(\bigcup_{l<k}[j_l,j_l+p]\right)\right)
 \ne\emptyset.
$$
Thus the lemma follows from the proof of \cite[Theorem~4.2]{A}
(with $n=p$ and $r=p$).
\endverif

The following lemma contains the main technical result needed to
prove Proposition~\ref{4.2}.

\begin{lemma} \label{4.5}
Let $N$ be a local algebra, $H\in N_{sa}$, $-\phi\in N^*$ is a tangent
functional to
$(P_\alpha)|_{N_{sa}}$ at $H$. Let $E\colon A\to N$ be a conditional
expectation. Then for any function $f\in\cal D$ (the space of
$C^\infty$-functions with compact support) and any $a,b\in N$ we have

$\displaystyle
\left|\int_\R\hat f(t)\phi(aE(\sigma^H_t(b)))dt
   -\int_\R\hat f(t+i)\phi(E(\sigma^H_t(b))a)dt\right|$
$$
   \le||a||\int_\R(|\hat f(t)|+|\hat f(t+i)|)
         ||\sigma^H_t(b)-E(\sigma^H_t(b))||dt.
$$
\end{lemma}

\noindent
{\it Proof.} First consider the case when $(P_\alpha)|_{N_{sa}}$ is
differentiable at $H$, in other words $-\phi$ is the unique tangent
functional. With the notations of Lemma~\ref{4.3} consider the state
$f_n$ on $\vee_{j\in I_n}\alpha^j(A_n)$ with density operator
$$
\left(\Tr_{\vee_{j\in I_n}\alpha^j(A_n)}
      \left(e^{-\sum_{j\in I_n}\alpha^j(H)}\right)\right)^{-1}
      e^{-\sum_{j\in I_n}\alpha^j(H)}.
$$
Then define a positive linear functional $\phi_n$ on $N$ by
$$
\phi_n(x)={1\over k_nm_n}\sum_{j\in I_n}f_n(\alpha^j(x)).
$$
Note that $||\phi_n||=\phi_n(1)\le1$. Since $-f_n$ is a tangent
functional to the convex function
$x\mapsto\log\Tr_{\vee_{j\in I_n}\alpha^j(A_n)}(e^{-x})$ on
$(\vee_{j\in I_n}\alpha^j(A_n))_{sa}$ at the point $\sum_{j\in
I_n}\alpha^j(H)$, $-\phi_n$ is a tangent functional to the function
$N_{sa}\ni x\mapsto{1\over k_nm_n}\log\Tr_{\vee_{j\in I_n}\alpha^j(A_n)}
      \left(e^{-\sum_{j\in I_n}\alpha^j(x)}\right)$ at $H$. It
follows that any limit point of the sequence $\{-\phi_n\}_n$ is a
tangent functional to $(P_\alpha)|_{N_{sa}}$ at $H$. Since the latter
is unique by assumption, $\phi_n\to\phi$ as $n\to\infty$.

Since $f_n$ is a $(\sigma^{H,I_n}_t,1)$-KMS state, by
\cite[Proposition~5.3.12]{BR2} we have
\begin{equation} \label{e4.2}
\int_\R\hat f(t)f_n(\alpha^j(a)\sigma^{H,I_n}_t(\alpha^j(b)))dt
   =\int_\R\hat
   f(t+i)f_n(\sigma^{H,I_n}_t(\alpha^j(b))\alpha^j(a))dt\ \
   \forall j\in I_n.
\end{equation}

Note that $\sigma^{H,I_n}_t(\alpha^j(b))=\alpha^j(\sigma^{H,I_n-j}_t(b))$.
Fix $q\in\N$, and set $I_{n,q}=\cup^{m_n-1}_{j=0}[jk_n+q,jk_n+k_n-p_n-q]$.
By Lemma~\ref{4.4}, if $q$ is large enough then
$\sigma^{H,I_n-j}_t(b)$ is arbitrarily close to $\sigma^H_t(b)$ for
any $j\in I_{n,q}$ and any $t$ in a fixed compact subset of~$\R$. But
then $\sigma^{H,I_n}_t(\alpha^j(b))-\alpha^j(E(\sigma^H_t(b)))$ is
arbitrarily close to $\alpha^j(\sigma^H_t(b)-E(\sigma^H_t(b)))$. In
other words,

$\displaystyle
\left|\int_\R dt\,\hat f(t){1\over k_nm_n}\sum_{j\in I_{n,q}}
   f_n\left(\alpha^j(a)\sigma^{H,I_n}_t(\alpha^j(b))
      -\alpha^j(a)\alpha^j(E(\sigma^H_t(b)))\right)\right|$
$$
   \le||a||\int_\R|\hat f(t)|\,||\sigma^H_t(b)-E(\sigma^H_t(b))||dt+\eps(q)\ \
\forall n\in\N,
$$
where $\eps(q)\to0$ as $q\to\infty$. Since $\# I_{n,q}/\# I_n\to1$ as
$n\to\infty$, letting $n\to\infty$ we may replace averaging over the
set $I_{n,q}$ by averaging over $I_n$, and then obtain

$\displaystyle
\limsup_{n\to\infty}\left|{1\over k_nm_n}\sum_{j\in I_n}\int_\R\hat f(t)
   f_n\left(\alpha^j(a)\sigma^{H,I_n}_t(\alpha^j(b))\right)dt
      -\int_\R\hat f(t)\phi(aE(\sigma^H_t(b)))dt\right|$
$$
   \le||a||\int_\R|\hat f(t)|\,||\sigma^H_t(b)-E(\sigma^H_t(b))||dt.
$$
Since an analogous estimate holds for
$\int_\R\hat f(t+i)\phi(E(\sigma^H_t(b))a)dt$, we obtain the
conclusion of the lemma by virtue of (\ref{e4.2}).

\medskip
If $(P_\alpha)|_{N_{sa}}$ is not differentiable at $H$ then by
\cite[Theorem~1]{LR}, $\phi$ lies in the closed convex hull of those
$\tilde\phi$, for which there exists a sequence $\{H_n\}_n\subset
N_{sa}$ converging to $H$ such that $(P_\alpha)|_{N_{sa}}$ has a
unique tangent functional $-\phi_n$ at $H_n$ and
$\phi_n\to\tilde\phi$. Since for $\phi_n$ the lemma is already proved
(for $H_n$ instead of~$H$), using Lemma~\ref{4.4} we conclude that
the conclusion of the lemma is true for $\tilde\phi$. But then it is true
for any functional in the closed convex hull of the~$\tilde\phi$'s.
\endverif

\noindent
{\it Proof of Theorem~\ref{4.2}.} If $-\phi$ is a tangent
functional for $P_\alpha$ at $H$ then $-\phi|_N$ is a tangent
functional for $(P_\alpha)|_{N_{sa}}$ at $H$ for any local algebra $N$
containing $H$. Thus by Lemma~\ref{4.5} the equality
$$
\int_\R\hat f(t)\phi(a\sigma^H_t(b))dt
   =\int_\R\hat f(t+i)\phi(\sigma^H_t(b)a)dt
$$
holds for all $f\in\cal D$ and all local $a,b$, hence for all $a,b\in
A$. By \cite[Proposition~5.3.12]{BR2} this is equivalent to the
KMS-condition.
\endverif

\begin{remark} \label{4.6}
\rm Under the assumptions of Theorem~\ref{4.1}, if
$$
\phi\in\bigcap_{\eps>0}\overline{
   \{\psi\,|\,P_\alpha(H)<h_\psi(\alpha)-\psi(H)+\eps\}}
$$
(weak$^*$ closure), then $-\phi$ is a tangent functional for the
pressure at $H$, hence $\phi$ is a $(\sigma^H_t,1)$-KMS state. In
other words, any weak$^*$ limit point of a sequence on which the sup in the
variational principle is attained, is a $(\sigma^H_t,1)$-KMS state.
If $ht(\alpha)=+\infty$, this is of course false in general. Moreover, for
any $\alpha$-invariant state $\phi$ there exists a sequence
$\{\phi_n\}_n$ converging in norm to $\phi$
such that $h_{\phi_n}(\alpha)=+\infty$ for all $n$. Indeed, first
note that by taking infinite convex combinations of states of large
entropy we can find a state $\psi$ of infinite entropy. Then $\phi_n={1\over
n}\psi+{n-1\over n}\phi\tends{n\to\infty}\phi$ and
$h_{\phi_n}(\alpha)=+\infty$.
\end{remark}

\bigskip\bigskip

\section{Examples} \label{5}

First we consider a class of systems arising naturally from systems of
topological dynamics.

Let $\sigma$ be an expansive homeomorphism of a
zero-dimensional compact space $X$, $G$ the group of uniformly
finite-dimensional homeomorphisms of~$X$ in the sense of Krieger
\cite{K}. By definition, a homeomorphism $T$ belongs to $G$ if
$$
\lim_{|n|\to\infty}\sup_{x\in X}d(\sigma^nTx,\sigma^nx)=0,
$$
where $d$ is a metric defining the topology of~$X$. In other words, $G$
consists of those homeomorphisms $T$ of~$X$, for which there exists a
bound on the number of coordinates of any point that are changed
under the action of~$T$, when $(X,\sigma)$ is represented as a
subshift by means of some generator. Since the group $G$ is locally
finite, the orbit equivalence relation $\RR\subset X\times X$ has a
structure of AF-groupoid \cite{Re}. Consider the groupoid
C$^*$-algebra $A=C^*(\RR)$ and the automorphism $\alpha$ of~$A$
defined by $\alpha(f)=f\circ(\sigma\times\sigma)$. The algebra $C(X)$
is a subalgebra of~$A$, and there exists a unique conditional
expectation $E\colon A\to C(X)$. Let $C_0(X)$ be the $*$-subalgebra of
$C(X)$ spanned by characteristic functions of clopen sets, and
$C_0(X,\R)$ the subalgebra of~$C_0(X)$ consisting of real functions.
Every element $g\in G$ defines a canonical unitary $u_g\in A$ such
that $u_gfu^*_g=f\circ g^{-1}$ for $f\in C(X)$. The $*$-algebra
generated by $C_0(X)$ and $u_g$, $g\in G$, is our algebra $\A$ of
local operators.

For $H\in C_0(X,\R)$ consider the $1$-cocycle $c_H\in Z^1(\RR,\R)$,
$$
c_H(x,y)=\sum_{j\in\Z}(H(\sigma^jx)-H(\sigma^jy)).
$$
Recall \cite[Definition~3.15]{Re} that a measure $\mu$ on
$X=\RR^{(0)}$ satisfies the $(c_H,1)$-KMS condition if its modular
function is equal to $e^{-c_H}$. In other words,
$$
{dg_*\mu\over d\mu}(x)=e^{-c_H(g^{-1}x,x)}.
$$

\begin{proposition} \label{5.1}
Let $H\in C_0(X,\R)$. Then

{\rm(i)} Any measure $\mu$ on $X$ which is an equilibrium measure
at~$-H$ satisfies the
$(c_H,1)$-KMS condition. In particular, any measure of maximal
entropy is $G$-invariant.

{\rm(ii)} The mapping $\mu\mapsto\mu\circ E$ defines a one-to-one
correspondence between equilibrium measures on $X$ at~$-H$ and
equilibrium states on $C^*(\RR)$ at $H$.
\end{proposition}

\noindent
{\it Proof.} First note that if $\phi$ is an $\alpha$-invariant state
on $C^*(\RR)$, and $\mu=\phi|_{C(X)}$ then
$$
h_\mu(\sigma)=h_{\mu\circ E}(\alpha)\ge h_\phi(\alpha).
$$
The equality is proved by standard arguments
using~\cite[Corollary~VIII.8]{CNT}.
The inequality follows from the fact that if $\psi$ is a state on a
finite dimensional C$^*$-algebra~$M$ with a masa~$B$ then $S(\psi)\le
S(\psi|_B)$. It follows that if $\mu$ is an equilibrium measure at
$-H$ then $\mu\circ E$ is an equilibrium state at $H$, and if $\phi$ is
an equilibrium state at $H$ then $\phi|_{C(X)}$ is an equilibrium
measure at~$-H$. By Theorem~\ref{4.1} any equilibrium state is
a $(\sigma^H_t,1)$-KMS state. But by \cite[Proposition~5.4]{Re} any
$(\sigma^H_t,1)$-KMS state has the form $\mu\circ E$ for some
measure $\mu$ satisfying the $(c_H,1)$-KMS condition. From this both
assertions of the proposition follow.
\endverif

\begin{example} \label{5.1a}
\rm As an application of Proposition~\ref{5.1} consider a topological
Markov chain $(X,\sigma)$  with transition matrix $A_T$. As is well-known,
if $A_T$ is primitive then the Perron-Frobenius theorem implies the
uniqueness of the trace on $C^*(\RR)$. If $A_T$ is only supposed to be
irreducible, then the traces of $C^*(\RR)$ form a simplex with the
number of vertices equal to the index of cyclicity of the matrix. The
barycenter of this simplex is the unique $\alpha$-invariant trace. By
Proposition~\ref{5.1} we conclude that if $A_T$ is irreducible then
$(X,\sigma)$ has a unique measure of maximal entropy. Thus we have
recovered a well-known result of Parry (see~\cite[Theorem~8.10]{W}).
\end{example}

While in the abelian case uniquely ergodic systems are of great
interest, they are not so for asymptotically abelian systems with locality.
Indeed, we have

\begin{proposition} \label{5.2}
Let $(A,\alpha)$ be a C$^*$-dynamical system which is asymptotically
abelian with locality. If there is a unique invariant state $\tau$,
then $\tau$ is a trace, and $\pi_\tau(A)$ is an abelian algebra.
\end{proposition}

For later use the main part of the proof will be given in a separate
lemma.

\begin{lemma} \label{5.3}
Let $(A,\alpha)$ be an asymptotically abelian system with locality,
$\tau$ an $\alpha$-invariant ergodic trace on $A$, $H$ a local
self-adjoint operator. Suppose for each $H'$ in the real linear span
of~$\alpha^j(H)$, $j\in\Z$, and for each $k\in\N$ there exists an
$\alpha^k$-invariant $(\sigma^{H',k\Z}_t,1)$-KMS state $\phi$ such that
$\tau={1\over k}\sum^{k-1}_{j=0}\phi\circ\alpha^j$. Then $\pi_\tau(H)$
is central in~$\pi_\tau(A)$.
\end{lemma}

\noindent
{\it Proof.} Replacing $A$ by $A/\Ker\,\pi_\tau$ we may identify $A$
with $\pi_\tau(A)\subset B({\cal H}_\tau)$.

The automorphism $\alpha$ being extended to $A''\subset B({\cal H}_\tau)$
is strongly asymptotically abelian. Hence for any $k\in\N$ the fixed
point algebra $(A'')^{\alpha^k}$ is central. Since $\alpha$ is
ergodic, this algebra is $k_0$-dimensional for some $k_0|k$, and we
may enumerate its atoms $z_1,\ldots,z_{k_0}$ in such a way that
$\alpha(z_1)=z_2,\ldots,\alpha(z_{k_0-1})=z_{k_0},\alpha(z_{k_0})=z_1$.
Now if $\phi$ is an $\alpha^k$-invariant state such that
$\tau={1\over k}\sum^{k-1}_{j=0}\phi\circ\alpha^j$ then $\phi\le
k\tau$, hence $\phi(x)=\tau(xa)$ for some positive
$a\in(A'')^{\alpha^k}$. In particular, $\phi$ is a trace. So if in
addition $\phi$ is a $(\sigma^{H',k\Z}_t,1)$-KMS state then the dynamics
$\sigma^{H',k\Z}_t$ is trivial on $(A'')_{s(a)}$, where $s(a)$ is
the support of~$a$. Hence $\delta_{H',k\Z}(y)s(a)=0$ for all local
$y$. Then $\delta_{H',k\Z}(y)z_i=0$ for some $z_i$ ($1\le i\le k_0$)
majorized by~$s(a)$.

Fix a local $x\in A$. Choose $p\in\N$ such that $\alpha^j(H)$ commutes
with $x$ whenever $|j|\ge p$. Pick any $m>p$ and set $k=2m+1$. For
$\lambda\in\R^k$ consider the operator
$$
H(\lambda)=\sum^m_{j=-m}\lambda_j\alpha^j(H).
$$
Applying the result of the previous paragraph to $H'=H(\lambda)$, we
find $i$, $1\le i\le k_0$, such that
$\delta_{H(\lambda),k\Z}(y)z_i=0$ for all local $y$. Denote by $X_i$
the set of all $\lambda\in\R^k$ satisfying the latter condition. Since
$\R^k=\cup_iX_i$, there exists $i$ for which $\R^k$ coincides with
the linear span of~$X_i$. Without loss of generality we may suppose
that $i=1$. Since for any $j\in[-m+p,m-p]$, any $j'\ne0$, and any
$\lambda\in X_1$, the elements $\alpha^{j'k}(H(\lambda))$ and
$\alpha^j(x)$ commute, we obtain
$$
0=\delta_{H(\lambda),k\Z}(\alpha^j(x))z_1=[H(\lambda),\alpha^j(x)]z_1,
$$
hence $[\alpha^{j'}(H),\alpha^j(x)]z_1=0$ for $j'\in[-m,m]$. In
particular,
$$
\alpha^j([H,x])z_1=0\ \ \hbox{for}\ \ j\in[-m+p,m-p].
$$
If $k_0\ne k$ and $k-2p\ge{k\over2}(>k_0)$ then
$\vee^{m-p}_{j=-m+p}\alpha^j(z_1)=1$, so $[H,x]=0$. If $k_0=k$
then $[H,x]z=0$, where
$$
z=\bigvee^{m-p}_{j=-m+p}\alpha^j(z_1)=\sum^{m-p}_{j=-m+p}\alpha^j(z_1), \ \
\tau(z)={k-2p\over k}.
$$
Since ${k-2p\over k}\to1$ as $m\to\infty$, we conclude that
$[H,x]=0$.
\endverif

\noindent
{\it Proof of Proposition~\ref{5.2}.} Since $A$ is a unital
AF-algebra, there exists a trace on $A$, hence there exists an
$\alpha$-invariant trace. It follows that the unique
$\alpha$-invariant state is a trace.

If $H$ is local then for any subset $I$ of~$Z$ there exists a
$(\sigma^{H,I}_t,1)$-KMS state. Indeed, if we take an increasing
sequence of finite subsets $I_n$ of~$I$ such that $\cup_nI_n=I$, an
increasing sequence of local algebras $A_n$ such that $\alpha^j(H)\in
A_n$ for $j\in I_n$ and $\cup_nA_n$ is dense in~$A$, and a sequence of
states $\phi_n$ such that $\phi_n|_{A_n}$ is a $(\sigma^{H,I_n}_t,1)$-KMS
state, then any weak$^*$ limit point of the sequence $\{\phi_n\}_n$ will be a
$(\sigma^{H,I}_t,1)$-KMS state. If in addition $I+k=I$, then the
state can be chosen to be $\alpha^k$-invariant (since the set of
$(\sigma^{H,I}_t,1)$-KMS states is $\alpha^k$-invariant). But if
$\phi$ is an $\alpha^k$-invariant state then the state ${1\over
k}\sum^{k-1}_{j=0}\phi\circ\alpha^j$ is $\alpha$-invariant, hence it
coincides with $\tau$. Thus the conditions of Lemma~\ref{5.3} are
satisfied. Hence $\pi_\tau(H)$ is central in~$\pi_\tau(A)$ for any
local $H$, so $\pi_\tau(A)$ is abelian.
\endverif

We consider two examples illustrating Proposition~\ref{5.2}.

\begin{example} \label{5.4}
\rm Let $U$ be the bilateral shift on a separable Hilbert space $\cal H$, and
$\alpha=\Ad U|_A$, where $A$ is the C$^*$-algebra $K({\cal H})+\C1$,
$K({\cal H})$ being the algebra of compact operators. Then the only
$\alpha$-invariant state is the trace $\tau$, which annihilates $K({\cal H})$.
Then $\pi_\tau(A)=\C1$.
\end{example}

\begin{example} \label{5.5}
\rm More generally, consider a uniquely ergodic system $(X,\sigma)$
and construct a system $(C^*(\RR),\alpha)$ as above. Let $\tau$ be
an $\alpha$-invariant trace. Then $\tau=\mu\circ E$ for some measure
$\mu$, and the unique ergodicity of~$(X,\sigma)$ means that $\mu$ is
the unique invariant measure. We check the conditions of Lemma
\ref{5.3} for any $H\in C(X_0,\R)$.

By the same reasons as in the proof of
Lemma~\ref{5.3}, the fixed point algebra $(\pi_\tau(A)'')^\alpha$ is
central. By \cite{FM} the center of the algebra $\pi_\tau(A)''$ is
isomorphic to $L^\infty(X,\mu)$. Since the measure $\mu$ is ergodic,
we conclude that the trace $\tau$ is also ergodic.

Let $H\in C(X_0,\R)$, $k\in\N$, and $\phi$ any $\alpha^k$-invariant
$(\sigma^{H,kZ}_t,1)$-KMS state. Then $\phi=\nu\circ E$ for some
$\sigma^k$-invariant measure $\nu$. Since
$\mu={1\over k}\sum^{k-1}_{j=0}\nu\circ\sigma^j$, we have
$\tau={1\over k}\sum^{k-1}_{j=0}\phi\circ\alpha^j$.

Thus we can apply Lemma~\ref{5.3}, and conclude that
$\pi_\tau(C(X))\cong C(\supp\mu)$ is central in~$\pi_\tau(A)$. This
means that $G$ acts trivially on $\supp\mu$, and $\pi_\tau(A)=C(\supp\mu)$.
By \cite[Proposition~4.5]{Re} the kernel of~$\pi_\tau$ is the algebra
corresponding to the groupoid $\RR_{X\backslash\supp\mu}$. Since there is
no non-zero finite $\sigma$-invariant measures on $X\backslash\supp\mu$,
any $\alpha$-invariant state is zero on $\Ker\,\pi_\tau$. Thus the
system $(A,\alpha)$ is uniquely ergodic and $\pi_\tau(A)=C(\supp\mu)$.
\end{example}

We next give an example of an asymptotically abelian C$^*$-dynamical
system $(A,\alpha)$ with $A$ an AF-algebra, for which there exist
non-tracial $\alpha$-invariant states with maximal finite entropy.
Hence the assumption of locality in Theorem~\ref{4.1} is essential.

\begin{example} \label{5.6}
\rm Let $\cal H$ be an infinite-dimensional Hilbert space, $A$ the
even CAR-algebra over~$\cal H$, $\alpha$ the Bogoliubov automorphism
corresponding to a unitary $U$. It easy to see that $\alpha$ is
asymptotically abelian if and only if $(U^nf,g)\tends{n\to\infty}0$ for any
$f,g\in\cal H$. If in addition $U$ has singular spectrum then by the
proof of \cite[Theorem~5.2]{SV} we have $ht(\alpha)=0$, while there
are many non-tracial $\alpha$-invariant states (for example,
quasi-free states corresponding to scalars $\lambda\in(0,1/2)$).
Unitaries with such properties can be obtained using Riesz products.
We shall briefly recall the construction.

Let $q>3$ be a real number, $\{n_k\}^\infty_{k=1}$ a sequence of
positive integers such that ${n_{k+1}\over n_k}\ge q$,
$\{a_k\}^\infty_{k=1}$ a sequence of real numbers such that
$a_k\in(-1,1)$, $a_k\to0$ as $k\to\infty$, $\sum_ka^2_k=\infty$. Then
the sequence of measures
$$
{1\over2\pi}\left[\prod^n_{k=1}(1+a_k\cos n_kt)\right]dt
$$
on $[0,2\pi]$ converges weakly$^*$ to a probability measure $\mu$ with
Fourier coefficients
$$
\hat\mu(n)=\mu(e^{int})=\cases{\displaystyle
   \prod^\infty_{k=1}\left({a_k\over2}\right)^{|\eps_k|},\ \
      \hbox{if} \ \ n=\sum_k\eps_kn_k \ \ \hbox{with} \ \
      \eps_k\in\{-1,0,1\},\cr
   \mbox{\ }\cr
   0,\ \ \hbox{otherwise}.}
$$
The measure $\mu$ is singular by \cite[Theorem~V.7.6]{Z}. We see also
that $\hat\mu(n)\to0$ as $|n|\to\infty$. Thus the operator $U$ of
multiplication by $e^{int}$ on $L^2([0,2\pi],d\mu)$ has the desired
properties.
\end{example}

\end{document}